\documentclass[11pt]{amsart}
\usepackage{syntonly}
\usepackage{hyperref}
\newcommand{\bbA}{{\mathbb A}}
\newcommand{\bbC}{{\mathbb C}}
\newcommand{\bbZ}{{\mathbb Z}}
\newcommand{\bbP}{{\mathbb P}}
\newcommand{\bbQ}{{\mathbb Q}}
\newcommand{\calV}{{\mathcal V}}
\newcommand{\calS}{{\mathcal S}}
\newcommand{\calC}{{\mathcal C}}
\newcommand{\calD}{{\mathcal D}}
\DeclareMathOperator{\Spec}{Spec}
\DeclareMathOperator{\HH}{H}
\DeclareMathOperator{\RR}{R}
\DeclareMathOperator{\CH}{CH}
\DeclareMathOperator{\Gal}{Gal}
\DeclareMathOperator{\tr}{tr}
\DeclareMathOperator{\et}{\text{\'et}}
\DeclareMathOperator{\Sym}{Sym}
\DeclareMathOperator{\IIJ}{IJ}
\DeclareMathOperator{\Ext}{Ext}
\DeclareMathOperator{\Div}{div}
\DeclareMathOperator{\cl}{cl}
\DeclareMathOperator{\GL}{GL}
\newcommand{\tensor}{\otimes}

\theoremstyle{plain}
\newtheorem{lemma}{Lemma}
\newtheorem{conj}{Conjecture}
\newtheorem{thm}[lemma]{Theorem}
\newtheorem{propn}[lemma]{Proposition}
\newtheorem{cor}[lemma]{Corollary}
\theoremstyle{definition}
\newtheorem*{ack*}{Acknowledgments}
\begin{document}

\title{Cycles over Fields of transcendence degree One}

\author[M. Green]{Mark Green}
\address{Dept. of Math., UCLA, Los Angeles, California, USA.}
\email{mlg@ipam.ucla.edu}

\author[P. A. Griffiths]{Philip A. Griffiths}
\address{IAS, Princeton, New Jersey, USA.}

\author[K. H. Paranjape]{Kapil H.~Paranjape}
\address{IMSc, CIT Campus, Tharamani, Chennai 600 113, India.}
\email{kapil@imsc.res.in}

\begin{abstract}
 We extend earlier examples provided by Schoen, Nori and Bloch to show
 that when a surface has the property that the kernel of its Albanese map
 is non-zero over the field of complex numbers, this kernel is non-zero
 over a field of transcendence degree one. This says that the conjecture
 of Bloch and Beilinson that this kernel is zero for varieties over
 number fields is precise in the sense that it is not valid for fields
 of transcendence degree one.
\end{abstract}

\maketitle

\section*{Introduction}

We work over subfields $k$ of $\bbC$, the field of complex
numbers.  For a smooth variety $V$ over $k$, the Chow group of
cycles of codimension $p$ is defined (see \cite{Fulton}) as
\[ \CH^{p}(V)=\frac{Z^{p}(V)}{R^{p}(V)}\]
where the group of cycles $Z^{p}(V)$ is the free abelian group on
scheme-theoretic points of $V$ of codimension $p$ and rational
equivalence $R^{p}(V)$ is the subgroup generated by cycles of the form
$\Div_{W}(f)$ where $W$ is a subvariety of $V$ of codimension $(p-1)$
and $f$ is a non-zero rational function on it. There is a natural
cycle class map
\[ \cl_{p}:\CH^{p}(V)\to \HH^{2p}(V)\]
where the latter denotes the singular cohomology group
$\HH^{2p}(V(\bbC),\bbZ)$ with the (mixed) Hodge structure given by
Deligne (see \cite{HodgeII}). The kernel of $\cl_{p}$ is denoted by
$F^{1}\CH^{p}(V)$. There is an Abel-Jacobi map (see
\cite{Griffiths}),
\[ \Phi_{p}:F^{1}\CH^{p}(V)\to \IIJ^{p}(\HH^{2p-1}(V))\]
where the latter is the intermediate Jacobian of a Hodge
structure, defined as follows
\[ \IIJ^{p}(H)=\frac{H\tensor \bbC}{F^{p}(H\tensor \bbC)+H}.\]
The kernel of $\Phi_{p}$ is denoted by $F^{2}\CH^{p}(V)$.
\begin{conj}[Bloch-Beilinson]
 If $V$ is a variety defined over a number field $k$ then
 $F^2\CH^p(V)=0$.
\end{conj}
We (of course) offer no proof of this conjecture. However, there are
examples due to C.~Schoen and M.~V.~Nori (see \cite{Schoen}), which
show that one cannot relax the conditions in this conjecture. In this
paper we present these and other examples to show that $F^2\CH^2(V)$
is non-zero for $V$ a variety over a field of transcendence degree at
least one, whenever it is non-zero over some larger algebraically
closed field.

We begin in Section~1 with a lemma. We then apply this lemma to prove
the following result by Hodge-theoretic methods for the second
Symmetric power of a curve.
\begin{propn}\label{MainProp}
  Let $C$ be a curve of genus at least two over a number field $F$.
  There is a non-torsion cycle $\xi$ in $F^2\CH^2(\Sym^2(C)\tensor
  K)$, where $K$ is the algebraic closure of the function field
  $F(C)$.
\end{propn}
In Section~2 we generalise this to surfaces other than $\Sym^2(C)$. To
do this with we need to use an $l$-adic analogue of the Abel-Jacobi
map. As before, let $V$ be a smooth projective variety over $k$. Let
$R$ be a finitely generated subring of $k$ for which there is a smooth
projective morphism $\pi:\calV\to\Spec(R)$, whose base change to $k$
gives $V$. Let $F^1\CH^p(\calV)$ denote the kernel of the natural
homomorphism
\[ 
   \CH^p(\calV)\to\HH^0_{\et}(\Spec(R),\RR^{2p}\pi_*\bbQ_l(p))
\]
We get a natural homomorphism
\[
  \Phi_{p,l,\calV}:F^{1}\CH^{p}(\calV)\to
       \HH^1_{\et}(\Spec(R),\RR^{2p-1}\pi_*\bbQ_l(p))
\]
Since any cycle (or
rational equivalence) on $V$ over $k$ is defined over a finitely
generated ring $R$, the direct limit of $F^1\CH^{p}(\calV)$ is
$F^1\CH^{p}(V)$. For want of better notation let us define
\[ 
 \IIJ^p_{\et}(\HH^{2p-1}(V)) = \lim_{\overrightarrow{R}} 
       \HH^1_{\et}(\Spec(R),\RR^{2p-1}\pi_*\bbQ_l(p))
\]
i.~e.\ $\IIJ^p$ is the corresponding direct limit of the
$\HH^1_{\et}$'s. An $l$-adic analogue of the Abel-Jacobi map is the
direct limit of the $\Phi_{p,l,\calV}$'s as we vary $R$,
\[
  \Phi_{p,l,\lim}:F^{1}\CH^{p}(V)\to \IIJ^p_{\et}(\HH^{2p-1}(V))
\]
On the other hand, it is clear that we have natural homomorphisms 
\[
     \HH^1_{\et}(\Spec(R),\RR^{2p-1}\pi_*\bbQ_l(p)) \to
       \Ext^1_{\Gal(\overline{k}/k)}(\bbQ_l(-p),\HH^{2p-1}_{\et}(V,\bbQ_l))
\]
One can show that these homomorphisms are injective when $k$ is the
field of fractions of the ring $R$ by finding extensions of $l$-adic
sheaves that represent elements of $\HH^1_{\et}$. The limiting
homomorphism
\[
  \Phi_{p,l}:F^{1}\CH^{p}(V)\to 
       \Ext^1_{\Gal(\overline{k}/k)}(\bbQ_l(-p),\HH^{2p-1}_{\et}(V,\bbQ_l))
\]
is thus another variant of the Abel-Jacobi homomorphism that has the
same kernel as $\Phi_{p,l,\lim}$ when $k$ is a finitely generated
field. We use the former map to define $F^2\CH^2(V)=\cap_l
\ker\Phi_{p,l,\lim}$. Recent results of Nori (unpublished) show that
the Hodge-theoretic definition of $F^2$ gives a group that is not
smaller than that given by the $l$-adic one above. It is still not
clear that the two definitions coincide as has been
conjectured by many (see \cite{Jannsen}), and has been proved by
Raskind~\cite{Wayne} for {\em zero cycles} by using the identification
of the Intermediate Jacobian with the Albanese variety in this case.

The following theorem uses the $l$-adic version of the definition of
$F^2$. In the rest of the paper the two definitions of $F^2$ will be
used without a change of notation; it will be clear from the context
which definition is being employed. By appealing to the result of
Raskind~\cite{Wayne} the two definitions coincide in our case.
\begin{thm}\label{MainThm}
  Let $S$ be any surface with $p_g(S)\neq 0$ defined over a number
  field $F$. There are non-trivial cycles in $F^2\CH^2(S\tensor
  K,\bbQ)$ where $K$ is the algebraic closure of the function field
  $F(T)$ in one variable $T$ over $F$.
\end{thm}

\begin{ack*}
Kapil Paranjape would like to thank Andreas Rosenschon for discussions
in which this problem was first raised. He would also like to thank
M. V. Nori and S. Bloch for discussions regarding the $\ell$-adic case. 
We thank W. Raskind for his careful reading of the paper and numerous
suggestions and corrections, especially regarding the use of limits in
\'etale cohomology.
\end{ack*}

\section{The Second Symmetric power of a Curve}
For a surface $S$, let $H^2(S)_{\tr}$ denote the orthogonal complement
in $H^2(S)$ of the (N\'eron-Severi) subspace generated by cohomology
classes of curves in $S$ which are defined over an algebraically
closed field.
\begin{lemma}[Green-Griffiths]\label{GG}
  Let $C$ be a smooth projective curve and $S$ be a smooth projective
  surface. Suppose $Z$ is a homologically trivial cycle of codimension
  two in $C\times S$. Assume that $C$, $S$ and $Z$ are defined over a
  subfield $k$ of $\bbC$. Suppose that {\em either one} of the following
  conditions hold:
  \begin{enumerate}
  \item The component in $\IIJ^2(\HH^1(C)\tensor\HH^2(S)_{\tr})$ of the
    Abel-Jacobi invariant of $Z$ is non-torsion. Or,
  \item The component in $\IIJ^2_{\et}(\HH^1(C)\tensor\HH^2(S)_{\tr})$
    of the $l$-adic Abel-Jacobi invariant $\Phi_{2,l,\lim}(Z)$ is
    non-torsion.
  \end{enumerate}
  Then the following map is non-torsion on $F^1\CH^1(C\tensor\bbC)$,
  \[ z=p_{2*}(p_1^*(\_)\cdot Z): \CH^1(C\tensor\bbC)\to
                        \CH^2(S\tensor\bbC)
  \]
\end{lemma}
\begin{proof}
  Let $A$ denote the canonical divisor of $C$ (which is defined over
  $k$). Let $g$ denote the generic point of $C$ considered as a point
  of $C\tensor k(C)$ and hence as a point of $C\tensor\bbC$ via some
  inclusion of $k(C)$ in $\bbC$. If $z(\deg(A)g-A)$ is torsion in
  $\CH^2(S\tensor\bbC)$, then a rational equivalence demonstrating
  this is defined over some finite extension of $k(C)$. Taking the
  norm of this equivalence we see that $z(\deg(A)g-A)$ is torsion in
  $\CH^2(S\tensor k(C))$ as well. Hence a multiple of $z(\deg(A)g-A)$
  is a sum of the form $\sum\Div_{\calD_i}(f_i)$ where
  ${\calD_i}\subset \Spec(k(C))\times S$ is a curve and $f_i$ is a
  function on ${\calD_i}$. Let $D_i\subset C\times S$ be the Zariski
  closure of $\calD_i$ and consider the cycle $W$ obtained by summing
  up $\Div_{D_i}(f_i)$. By assumption, $W$ restricts to a multiple of
  $z(\deg(A)q-A)$ on $\Spec(k(C))\times S$. Thus there are closed
  points $q_j$ of $C$ over $k$ and divisors $E_j$ in
  $\Spec(k(q_j))\times S$ so that for some integer $n$ we have
  \[
     \sum \Div_{D_i}(f_i) = W = n(\deg(A)Z - C\times z(A)) +
                            \sum_j q_j\times E_j
  \]
  The image of cycles of the form $C\times z(A)$ as well as cycles of
  the form $\sum_j q_j\times E_j$ in the relevant piece of $\IIJ^2$ is
  zero in case (1) as well as (2). This contradicts the hypothesis on
  $Z$.
\end{proof}  
From this lemma it follows that in order to construct a non-trivial
cycle in $F^2\CH^2(S)$ over a field of transcendence degree one it is
sufficient to construct $C$, $S$ and $Z$ as above over a field of
transcendence degree 0.
\begin{cor}
  Let $C$ be a smooth projective curve and $S$ a smooth projective
  surface over a subfield $F$ of $\overline{\bbQ}$. Let $Z$ be a
  homologically trivial cycle on $C\times S$ such that the induced
  homomorphism
  \[ z=p_{2*}(p_1^*(\_)\cdot Z): \CH^1(C\tensor\bbC)\to
                        \CH^2(S\tensor\bbC)
  \]
  is non-torsion. Then, there is a non-torsion class in
  $F^2\CH^2(S\tensor_F\overline{F(C)})$.
\end{cor}
\begin{proof}
  From a result of Roitman \cite{Roitman} it follows that the kernel
  of $z$ is a countable union of translates of proper abelian
  subvarieties $A_i$ of $J(C)$. But each such translated abelian subvariety is
  defined over $\overline{\bbQ}$ since $Z$ is defined over
  $\overline{\bbQ}$.  Using a point $p$ of $C$ over $\overline{\bbQ}$
  we embed $C$ in $J(C)$. Since the image of $C$ generates $J(C)$, the
  image of the generic point $\Spec(F(C))\to C$ is not in $A_i$ for
  any $i$. Its image under $z$ gives us the desired cycle in
  $F^2\CH^2(S\tensor_F\overline{F(C)})$.
\end{proof}
Let $C$ be a curve of genus at least 2 over $k=\overline{\bbQ}$. Let
$S=\Sym^2(C)$ and consider the cycle $X$ in $C\times C\times S$ which
is the graph of the natural map $q:C\times C\to S =\Sym^2(C)$. The map
$q_*$ on $\HH^{1,0}(C)\tensor\HH^{1,0}(C)$ identifies the range
$\HH^{2,0}(S)$ with the second exterior power of $\HH^{1,0}(C)$.  The
cohomology class $[X]$ therefore induces a non-zero map
\[
   p_{23*}(p_1^*(\_)\cup[X]): \HH^{1,0}(C) \to
                  \HH^{0,1}(C)\tensor\HH^{2,0}(S)
\]
Thus the map
\[
  J(C\tensor\bbC) = \IIJ^1(\HH^1(C)) \to
               \IIJ^2(\HH^1(C)\tensor\HH^2(S)_{\tr})
\]
is non-torsion. Hence, the image is a non-zero abelian variety $A$
which is a quotient of $J(C\tensor\bbC)$. Now, $J(C)$ is defined over
$\overline{\bbQ}$ and hence so is any quotient abelian variety. In
particular, $A$ is $B\tensor\bbC$ for some abelian variety $B$
defined over a number field.
\begin{propn}
  Let $J\to B$ be a surjective homomorphism of abelian varieties over
  a number field; $B\neq 0$. There are points in $J(\overline{\bbQ})$
  whose image in $B$ is non-torsion.
\end{propn}
\begin{proof}
  Let $L$ be the completion of the number field at a place lying over
  over a prime $p$ and $M$ be the algebraic closure of $\bbQ$ in $L$.
  $J(L)$ and $B(L)$ are $p$-adic Lie groups. Consider a finite
  morphism $\pi:J\to\bbP^n$ and let $U\subset\bbP^n$ be a Zariski open
  set over which $\pi$ is \'etale. It is clear that $U(M)$ is dense in
  $U(L)$. By the implicit function theorem the map
  $J(L)\cap\pi^{-1}(U)\to U(L)$ is open; moreover, the former set is
  dense in $J(L)$. Since $\pi$ is finite, any point in
  $J(L)\cap\pi^{-1}(U(M))$ is in $J(M)$. Thus $J(M)$ is dense in
  $J(L)$. Now $J(L)\to B(L)$ is open, so the image of $J(M)$ in
  $B(L)$ is dense in an open subgroup. The torsion subgroup of a
  $p$-adic Lie group is finite. The result follows.
\end{proof}
\begin{cor}
  Let $A$ be a non-zero abelian variety over a number fields. The
  group $A(\overline{\bbQ})$ is non-torsion.
\end{cor}
By this proposition and the discussion above we obtain a cycle $Z$ on
$C\times S$ over $\overline{\bbQ}$ as required. Now applying the lemma
we see that we have proved Proposition~\ref{MainProp}.

\section{The General Case}

We use an argument of Terasoma (as explained to us by Bloch and Nori)
to construct a triple $(C,S,Z)$ as in Lemma~\ref{GG} for {\em any}
surface $S$ with $p_g\neq 0$.
\begin{lemma}\label{Specl}
  Assume $(\tilde{C},\tilde{S},\tilde{Z})$ is a triple over some
  finitely generated subfield $k$ of $\bbC$ which satisfies the
  condition {\rm (2)} of Lemma~{\rm\ref{GG}}. Then, there exists a
  specialisation $(C,S,Z)$ over a number field
  $F\subset\overline{\bbQ}$ which satisfies the condition {\rm (2)} of
  Lemma~{\rm\ref{GG}}.
\end{lemma}
\begin{proof}
  Let $T$ be an affine variety over a number field contained in $k$
  and $(\calC,\calS,\Xi)$ be a triple over $S$ so chosen that this
  data restricts to $(\tilde{C},\tilde{S},\tilde{Z})$ under the
  natural map $\Spec k\to T$. Using the Gysin sequence, we obtain an
  exact sequence of $\bbQ_l$ constructible local systems on $T$,
  \[
     0 \to \RR^1_{\et}\pi_{T*}(\calC,\bbQ_l)\tensor
           \RR^2_{\et}\pi_{T*}(\calS,\bbQ_l)_{\tr}
           \to \calV \to \bbQ_l(-2)_T \to 0
  \]
  here $\RR^2_{\et}\pi_{T*}(\calS,\bbQ_l)_{\tr}$ denotes the local
  system on $T$ associated with the Galois representation
  $\HH^2_{\et}(\tilde{S}\tensor\overline{k},\bbQ_l)_{\tr}$. Moreover,
  the extension class of this short exact sequence is naturally
  identified with the component in
  \( 
   \HH^1_{\et}(T,
           \RR^1_{\et}\pi_{T*}(\calC,\bbQ_l)\tensor
           \RR^2_{\et}\pi_{T*}(\calS,\bbQ_l)_{\tr}\tensor
           \bbQ_l(2))
  \)
  of the class $\Phi_{2,l,\calC\times_T\calS}(\Xi)$. By assumption,
  condition (2) of Lemma~{\ref{GG}} means this component is non-zero. So
  the sequence is non-split.
  
  If $p$ is a closed point of $T$, let $C_p$, $S_p$, $Z_p$ denote the
  specialisations. By restricting the above sequence we obtain an
  exact sequence of representations of $\Gal(\overline{\bbQ}/k(p))$ or
  equivalently, of a decomposition group of $p$ in $\Gal(\overline{k}/k)$,
  \[
     0 \to \HH^1_{\et}(C_p\tensor\overline{\bbQ},\bbQ_l)\tensor
           (\RR^2_{\et}\pi_{T*}(\calS,\bbQ_l)_{\tr})_{\overline{k(p)}}
           \to \calV_{\overline{k(p)}} \to \bbQ_l(-2) \to 0
  \]
  There are two possible problems. 
  \begin{enumerate}
  \item There may be more divisors in $S_p$ than on the generic fibre
    of $\calS\to T$ so that
    \[
           \HH^2_{\et}(S_p\tensor{\overline{k(p)}},\bbQ_l)_{\tr}) 
           \text{~is a proper subset of~}
           (\RR^2_{\et}\pi_{T*}(\calS,\bbQ_l)_{\tr})_{\overline{k(p)}}
    \]
  \item The above exact sequence may split so that the extension class
   becomes zero for $Z_p$.
  \end{enumerate}
  These two problems are resolved by applying the following lemma to
  the $l$-adic representation
  \[ \calV_{\overline{k}}\oplus\HH^2_{\et}(S\tensor\overline{k},\bbQ_l) \]
\end{proof}
\begin{lemma}[Terasoma]
  Let $V$ be a continuous $l$-adic representation of the fundamental
  group of a variety $T$ over a number field. Then there are
  infinitely many closed points $p$ of $T$ so that the image of
  $\Gal(\overline{k(T)}/k(T))$ in $\GL(V)$ is the same as the image of
  the decomposition group of $p$.
\end{lemma}
\begin{proof}
  Let $G$ denote the image of $\Gal(\overline{k(T)}/k(T))$ in $\GL(V)$
  and $D_p$ denote the image of a decomposition group associated with
  a closed point $p$ of $T$. These are compact subgroups of $\GL(V)$.
  As in~\cite{Terasoma}, there is a subgroup $H$ of finite index in
  $G$ so that $D_p=G$ if and only if the natural map $D_p\to G/H$ is a
  surjection. Equivalently, if $f:T'\to T$ denotes the \'etale cover
  of $T'$ corresponding to this finite quotient of the fundamental
  group of $T$, then it is enough that $f^{-1}(p)$ is a closed point
  of $T'$. To find such a point $p$, we may shrink $T$ to assume that
  $T$ is affine and that there is an \'etale map $g:T\to\bbA^n$. Thus
  it is enough to find a point $q$ in $\bbA^n$ so that $(g\circ
  f)^{-1}(q)$ is a closed point of $T'$. Such a point can be found by
  the Hilbert irreducibility theorem.
\end{proof}

Now we consider a surface $S$ over $\overline{\bbQ}$ with $p_g(S)\neq
0$; let $p$ be a closed point on $S$. Murre (see \cite{Murre}) has
constructed a decomposition of the diagonal of a surface $S$,
\[ \Delta_S = p\times S + S\times p + X_{2,2} + X_{1,3} + X_{3,1} \]
for the diagonal $\Delta_S$ in $\CH^2(S\times S)\tensor\bbQ$; where
$X_{i,j}$ is a cycle class so that its cohomology class has a
non-zero K\"unneth component only in
$\HH^i(S)\tensor\HH^j(S)\tensor\bbQ$. Moreover, by the Hodge index
theorem we can write $nX_{2,2}=\sum_i C_i\times D_i + X_{2,2,\tr}$
where $C_i$ and $D_i$ are curves on $S$ and $X_{2,2,\tr}$ has its
cohomology class in $\HH^2(S)_{\tr}\tensor\HH^2(S)_{\tr}$.

Let $H\subset S\times\bbP^1$ be any pencil; $Y_{i,j}$ denote the cycles
on $H\times S$ that are the natural pull-backs of $X_{i,j}$ for each
$i$ and $j$; similarly, let $Y_{2,2,\tr}$ be the pull-back of
$X_{2,2,\tr}$. Let $H_0$ denote the divisor class of the fibres of the
pencil and $A$ the class of an ample divisor on $S$. The class
$Y_{2,2,\tr}$ on $H\times S$, restricts to a homologically trivial
class on every member $H_t$ of the pencil; however,
$\HH^2(H)_{\tr}=\HH^2(S)_{\tr}$ so that the cohomology class of
$Y_{2,2,\tr}$ in $\HH^2(H)_{\tr}\tensor\HH^2(S)_{\tr}$ is non-zero.

For any open subset $U\subset\bbP^1$, let $H_U=H\times_{\bbP^1} U$;
the map $\HH^2(H)_{\tr}\to\HH^2(H_U)$ is injective. Thus the pull-back
of $Y_{2,2,\tr}$ to $H_U\times S$ is not homologically trivial for all
$U$. By the Leray spectral sequence for the map $\pi:H_U\to U$, we see
that $\HH^2(H)_{\tr}$ maps injectively into
$\HH^1(U,\RR^1\pi_*(\bbQ_l))$.  Let $g$ denote the generic point of
$\bbP^1$. The restriction of $Y_{2,2,\tr}$ to $H_g\times S$ is a cycle
$Z_g$ such that the triple $(H_g,S\tensor_k k(g),Z_g)$ satisfies the
condition (2) of Lemma~\ref{Specl} over the field $k(g)$.
Applying Lemma~\ref{Specl} we see that we have a point $t$ in
$\bbP^1(\bbQ)$ so that $(H_t,S,Z_t)$ is a triple satisfying condition
3 of the Lemma~\ref{GG} of Green and Griffiths. Thus we have proved
Theorem~\ref{MainThm}.


\end{document}